\documentclass[12pt]{article}
\input tcilatex
\usepackage{amsfonts,epsf}

\newcommand{\psdiag}[3]{\hspace{1mm}\raisebox{-#1mm}{\epsfysize#2mm
\epsffile{#3.eps}}\hspace{1mm}}
\newcommand{\psbild}[2]{\vspace{3mm}\\ \mbox{\epsfysize=#1mm\epsffile{#2.eps}}
\vspace{3mm}}

\begin{document}

\author{Rui Pedro Carpentier\\
\\
{\small\it Departamento de Matem\'{a}tica and Centro de
Matem\'{a}tica Aplicada}\\
{\small\it  Instituto Superior T\'{e}cnico}\\
{\small\it Avenida Rovisco Pais, 1049-001 Lisboa}\\
{\small\it Portugal}}
\title{Topological notions for Kauffman and
Vogel's polynomial}
\date{25rd July, 2001}
\maketitle

\begin{abstract}

In \cite{2} Kauffman and Vogel constructed a rigid vertex regular
isotopy invariant for unoriented four-valent graphs embedded in
three dimensional space. It assigns to each embedded graph $G$ a
polynomial, denoted $[G]$,  in three variables, $A$, $B$ and $a$,
satisfying the skein relations: $$ [\psdiag{2}{6}{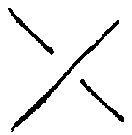}]=A
[\psdiag{2}{6}{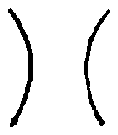}]+B [\psdiag{2}{6}{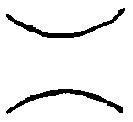}]+
[\psdiag{2}{6}{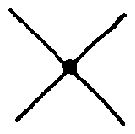}]$$ $$ [\psdiag{2}{6}{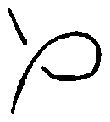}]=
a[\psdiag{2}{6}{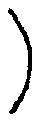}], \hspace{1cm}
[\psdiag{2}{6}{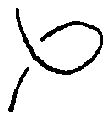}]= a^{-1} [\psdiag{2}{6}{straight}] $$
and is defined in terms of a state-sum and the Dubrovnik
polynomial for links.

In \cite{4} it is proved, in the case $B=A^{-1}$ and $a=A$, that
for a planar graph $G$ we have $[G]=2^{c-1}(-A-A^{-1})^v $, where
$c$ is the number of connected components of $G$ and $v$ is the
number of vertices of $G$.

In this paper we will show how we can calculate the polynomial,
with the variables $B=A^{-1}$ and $a=A$, without resorting to the
skein relation.

\end{abstract}

\section{Introduction.}

In \cite{1} a polynomial invariant is described for 4-valent rigid
vertex embedded graphs. This polynomial, denoted $[G]$ for a graph
$G$, with variables $A$, $B$ and $a$ satisfies the skein relation
$$[\psdiag{2}{6}{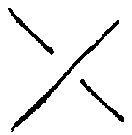} ]=A[\psdiag{2}{6}{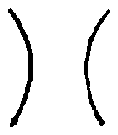} ]+B[\psdiag{2}{6}{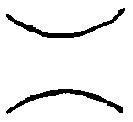}
]+[\psdiag{2}{6}{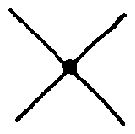} ]$$ and restricted to a link gives the
Dubrovnik polynomial with variables $a$ and $A-B$. It is proved in
\cite{2} that this polynomial is invariant for regular isotopies
and satisfies the conditions $$[\psdiag{2}{6}{overtwist}]=
a[\psdiag{2}{6}{straight}], [\psdiag{2}{6}{undertwist}]= a^{-1}
[\psdiag{2}{6}{straight}].$$

This means that the polynomials of graph diagrams are the same if
we can change one graph diagram to the other by the Reidemeister
moves for 4-valent graphs except for the first move (see
\cite{1,3}).

The non-invariance under the first Reidemeister move can be
corrected using the {\it twisting number} of a graph. This is
defined as the sum of the writhes over all knot-theoretic circuits
of the graph. A {\it knot-theoretic circuit} of a graph is a
closed walk on the graph which corresponds to a link component if
we change all vertices into crossings, and its {\it writhe} is the
sum of the signs of all its self-crossings. The twisting number of
a graph is invariant under regular isotopies. Thus, defining
$t(G)$ as the twisting number of a graph $G$, we have that
$a^{-t(G)}[G]$ is invariant under any isotopy since
$$t(\psdiag{2}{6}{overtwist})= t(\psdiag{2}{6}{straight})+1$$ and
$$t(\psdiag{2}{6}{undertwist})= t(\psdiag{2}{6}{straight})-1$$.

In \cite{2} the following graphical calculus is proved.

\begin{theorem}
For 4-valent graph diagrams, differing only in the marked local
picture, we have the following identities:

\begin{description}
\item  $[\psdiag{2}{6}{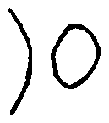} ]=\mu [\psdiag{2}{6}{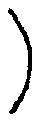}]$

\item  $[\psdiag{2}{6}{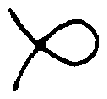} ]={\cal O}[\psdiag{2}{6}{g1} ]$

\item  $[\psdiag{2}{6}{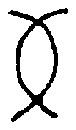} ]=(1-AB)[\psdiag{2}{6}{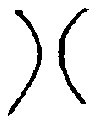} ]+\gamma
[\psdiag{2}{6}{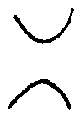} ]-(A+B)[\psdiag{2}{6}{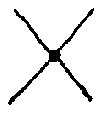} ]$

\item  $[\psdiag{2}{6}{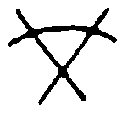}]-[\psdiag{2}{6}{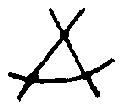}
]=AB([\psdiag{2}{6}{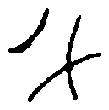}]-[\psdiag{2}{6}{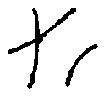}]+[\psdiag{2}{6}{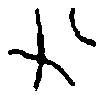}]-[\psdiag{2}{6}{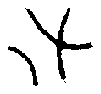}
]+[\psdiag{2}{6}{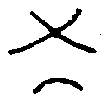}]-[\psdiag{2}{6}{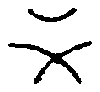}])+\xi
([\psdiag{2}{6}{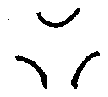}]-[\psdiag{2}{6}{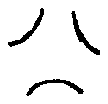}])$
\end{description}

where

\begin{description}
\item $\mu = \frac{a-a^{-1}}{A-B}+1$
\item ${\cal O}=\frac{Aa^{-1}-Ba}{A-B}-(A+B)$
\item $\gamma = \frac{B^2a-A^2a^{-1}}{A-B}+AB$
\item $\xi =\frac{B^3a-A^3a^{-1}}{A-B}$
\end{description}

\end{theorem}

As was shown in \cite{4}, this graphical calculus is sufficient to
calculate the polynomial of a planar graph.

In \cite{4} the author studied the case $B=A^{-1}$ and $a=A$. In
this case, the graphical
calculus takes the form:%

\begin{description}
\item $\lbrack \psdiag{2}{6}{k2} ]=2[\psdiag{2}{6}{g1} ]$,
\item $[\psdiag{2}{6}{l2} ]=-(A+A^{-1})[\psdiag{2}{6}{g1} ]$,
\item $\lbrack \psdiag{2}{6}{m2} ]=-(A+A^{-1})[\psdiag{2}{6}{l9} ]$,
\item $\lbrack
\psdiag{2}{6}{b3}]-[\psdiag{2}{6}{a3}]=[%
\psdiag{2}{6}{c3}]-[\psdiag{2}{6}{e3}]+[\psdiag{2}{6}{g3}]-[\psdiag{2}{6}{f3}]+
[\psdiag{2}{6}{h3}]-[\psdiag{2}{6}{i3}]-(A+A^{-1})([\psdiag{2}{6}{j3}]-[\psdiag{2}{6}{l4}])$,
\end{description}

and for planar graphs we have:

\begin{theorem}
In the case $a=A$ and $B=A^{-1}$ for any 4-valent planar graph $G$ we have $%
[G]=2^{c-1}(-A-A^{-1})^v$, where $c$ is the number of connected components
of $G$ and $v$ is the number of vertices of $G.$
\end{theorem}

This means that a necessary condition for a 4-valent rigid vertex
embedded graph $G$ to be isotopic to a planar graph is that
the polynomial of $G$ with $B=A^{-1}$ and $a=A$ is $%
[G]=2^{c-1}(-A-A^{-1})^vA^{t(G)}$ where $c$ is the number of
connected components, $v$ is the number of vertices and $t(G)$ is
the twisting number of $G$. This corollary is a consequence of the
fact that the polynomial $ a^{-t(G)}[G] $ is an isotopy invariant.

Using the skein relation $$[\psdiag{2}{6}{a1}
]=A[\psdiag{2}{6}{d1} ]+A^{-1}[\psdiag{2}{6}{e1}
]+[\psdiag{2}{6}{c1} ]$$ and the previous theorem it is easy to
prove the following result:

\begin{corollary}
If a graph diagram $G$ has only one crossing, and the removal of this crossing
does not
change the number of connected components, then the polynomial of $G$
vanishes in the case $B=A^{-1}$ and $a=A$, and thus $G$ is not isotopic to a
planar graph.
\end{corollary}

Thus we have a class of embedded graphs with polynomial equal to
zero and therefore not isotopic to planar graphs. In the next
section we will characterize the graphs with polynomial equal to
zero.

\section{Calculating the polynomial for embedded graphs.}

Let a {\it separating curve} be a simple regular curve, on the
plane, that intersects the diagram only at the crossings of the
diagram and does so in the following way: \psdiag{2}{6}{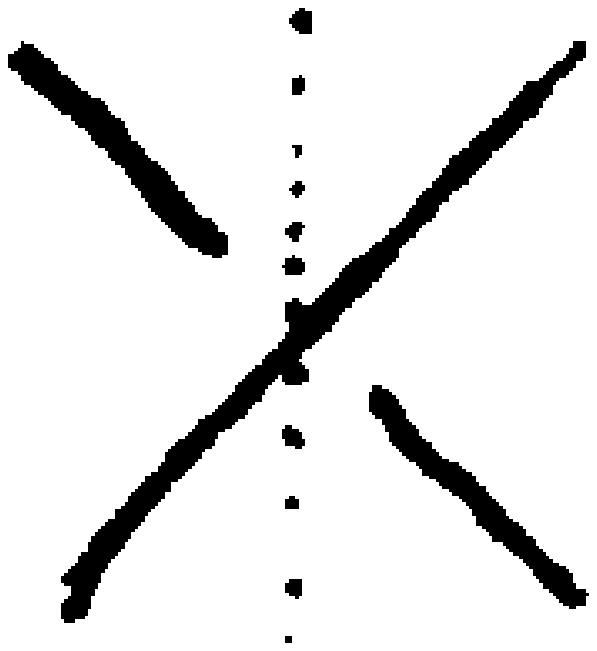} or
\psdiag{2}{6}{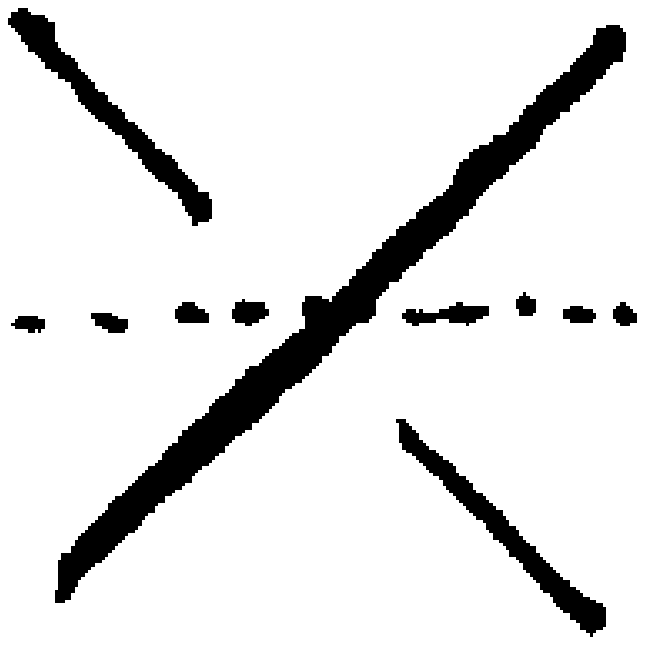}. A {\it partition} of the diagram is a choice of
marks of the type \psdiag{2}{6}{w5} or \psdiag{2}{6}{w6}, at each
crossing, such that there exists a collection of disjoint
separating curves intersecting all crossings of the diagram in
agreement with these marks. A partition without marks we call a
{\it null partition}\footnote{This only exists for diagrams
without crossings, i.e., planar graphs.}.

A graph diagram is called {\it separable} if it has a partition.
We can show that this property is invariant under ambient
isotopies. The planar graphs are clearly separable, since they
have the null partition. Later we will see that any link is also
separable.

\begin{proposition}
The polynomial, with variables $B=A^{-1}$ and $a=A$, of a graph
diagram is zero if and only if the graph diagram is not separable.
\end{proposition}

This result is a corollary of theorem 6.

There is a way to determine the partition of a graph. For that we
need to introduce another definition.

A orientation on the edges of a graph will be called {\it
hyperbolic} if at each vertex of the graph the orientation is of
the type \psdiag{2}{6}{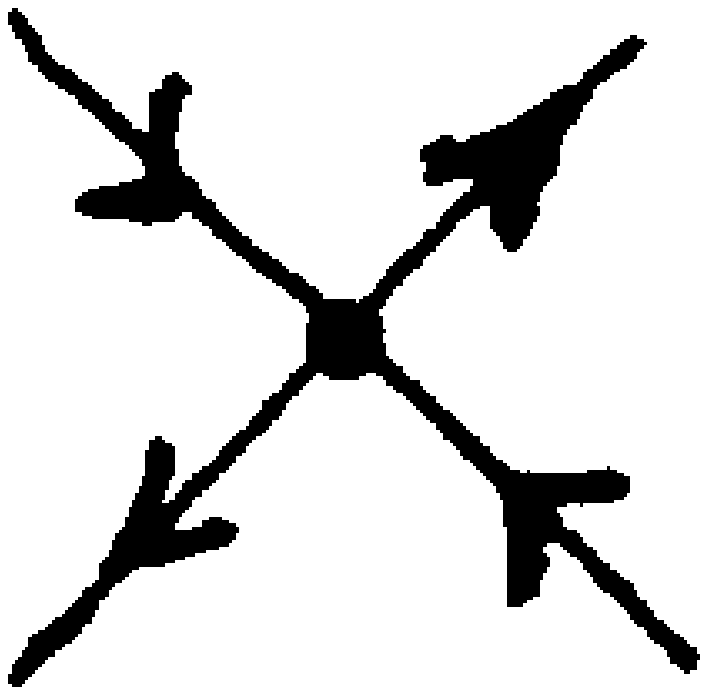}.

\begin{theorem}
For a given graph diagram, we have that a choice of marks on the
crossings is a partition if and only if there exists a hyperbolic
orientation on the edges of the graph such that the mark on each
crossing is in agreement with the orientation in the following way
\psdiag{2}{6}{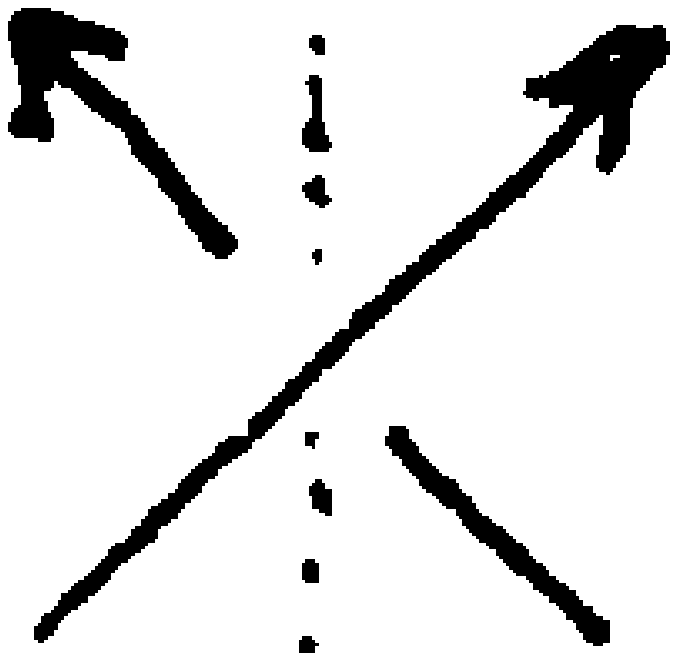} or \psdiag{2}{6}{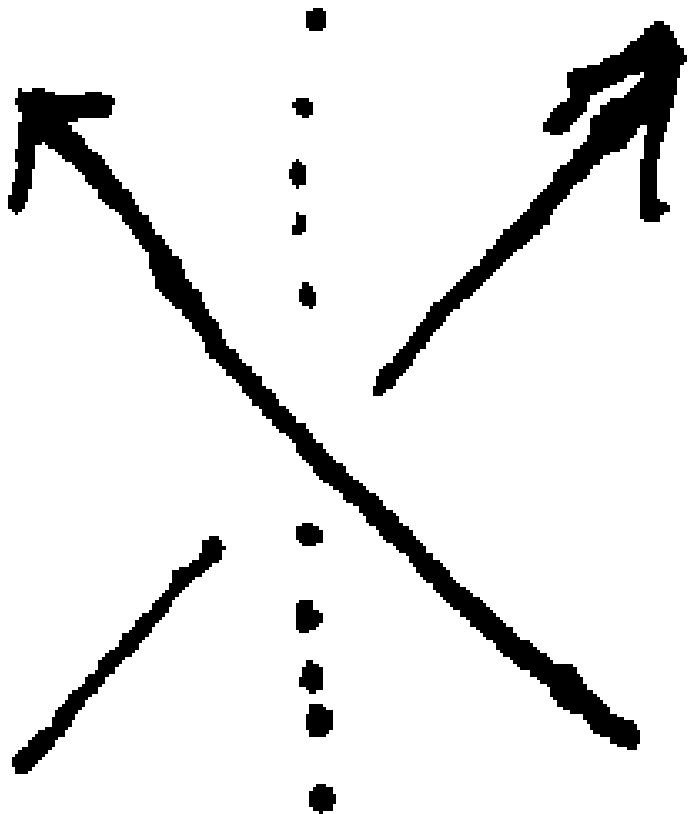}.

In particular, we have that a graph is separable if and only if it
has hyperbolic orientations.
\end{theorem}

\TeXButton{Proof}{\proof} Given a partition there exists a
collection of separating curves that cross the graph diagram in
agreement with the marks of the partition. These curves, together
with the diagram, form the skeleton of a map that is
bicolorable\footnote{i.e. the regions delimited by the diagram and
the curves can be colored with two colors (black and white) in
such way that two adjacent regions to the same edge of the diagram
or part of the curve have different colors (see \cite{5}).} since
the vertices of the skeleton of the map have even degree (four in
the case of the vertices of the diagram, and six in the case of
the crossings of the diagram cut by the curves). Thus we can
choose opposite orientations for the boundaries of the regions
with different colors. This gives an orientation on the diagram
that is constant on each edge (even when the edge through a
crossing), is hyperbolic and such that at each crossing the
separating curve cuts the diagram in the following way
\psdiag{2}{6}{w7} or \psdiag{2}{6}{z8}.

\centerline{\psbild{25}{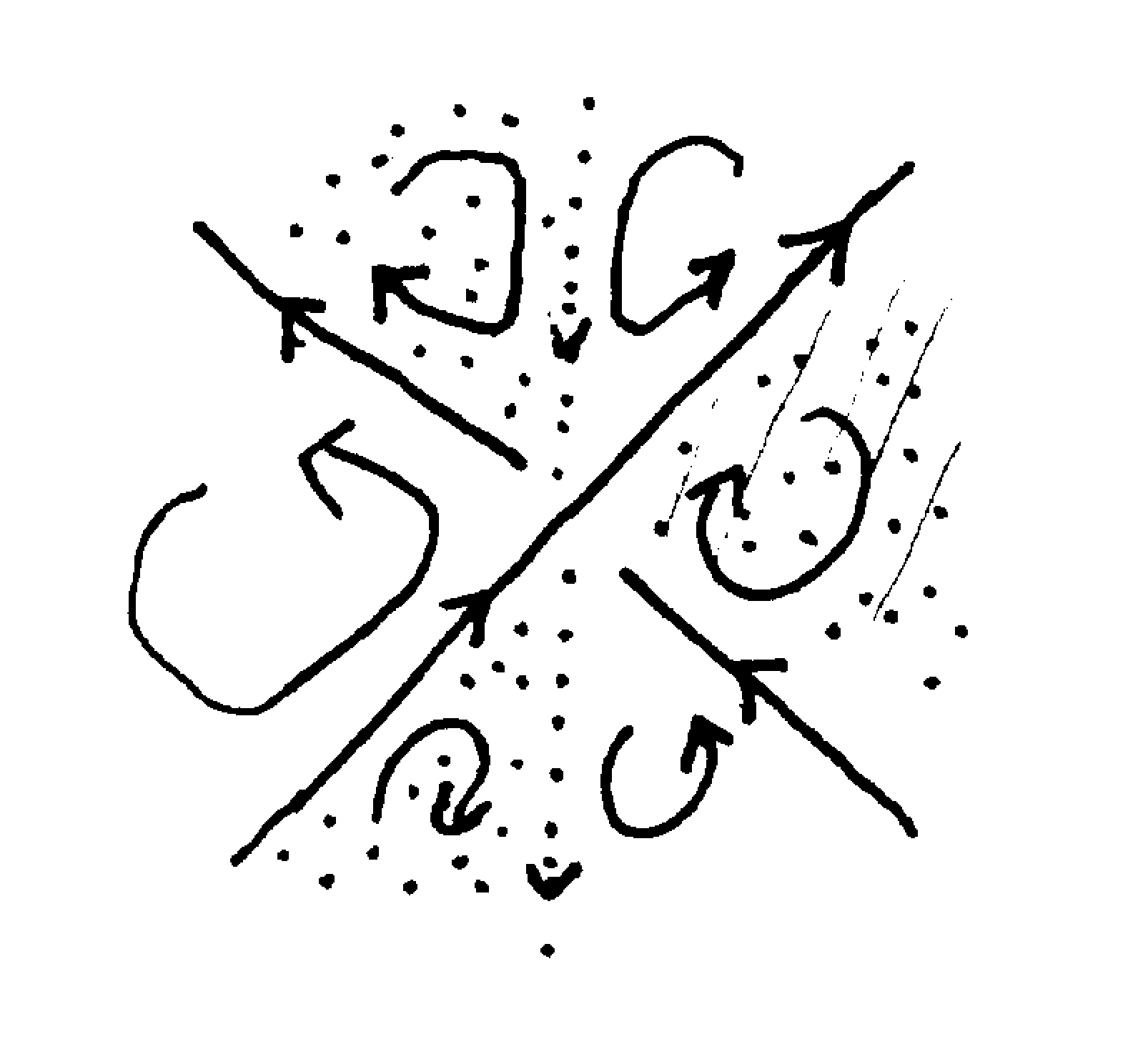}, \psbild{25}{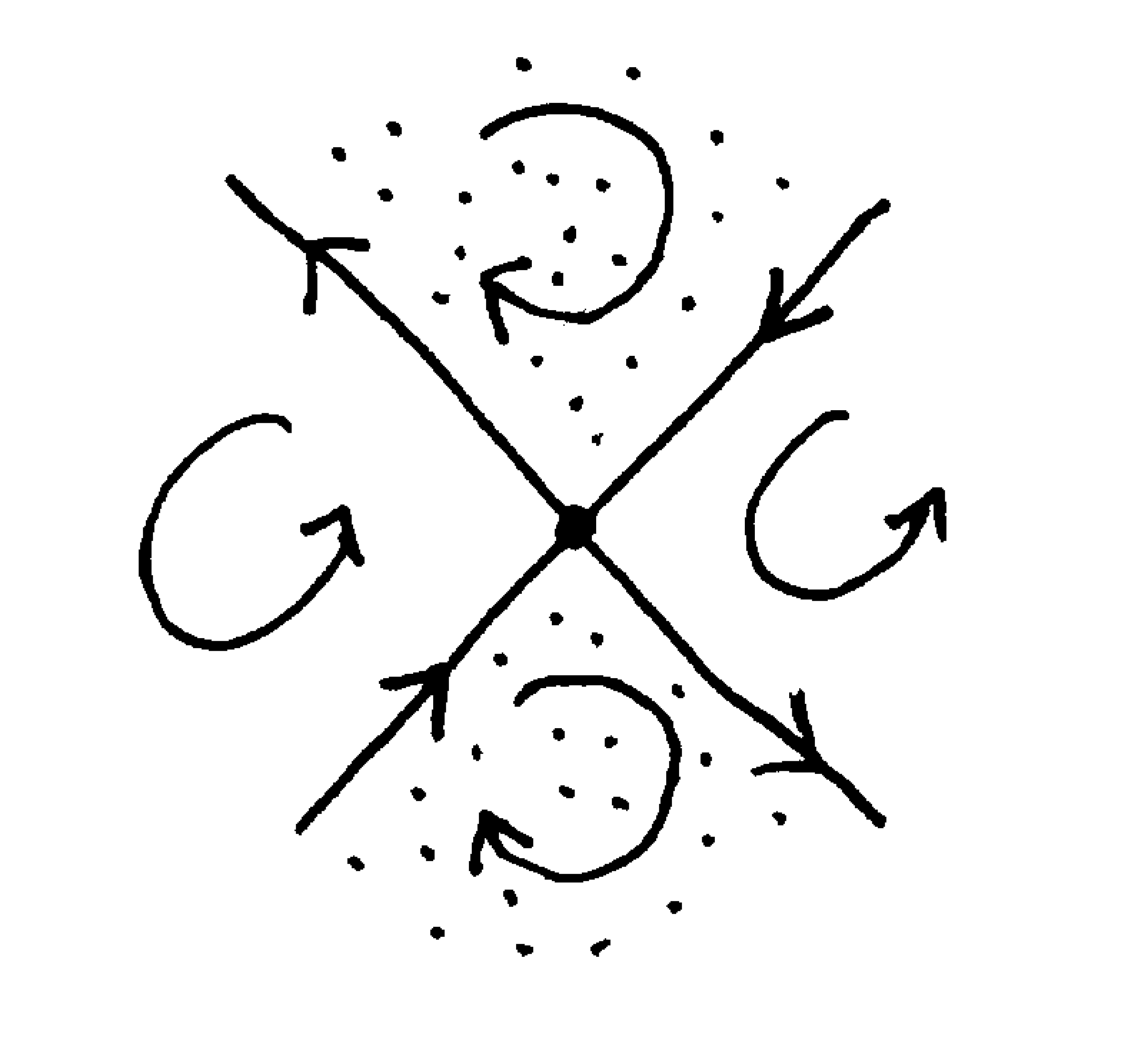}}

Given a hyperbolic orientation on the edges of the graph, we want
to see that the marks \psdiag{2}{6}{w7} or \psdiag{2}{6}{z8} form
a partition. For that we take the regions delimited by the
diagram. For each region the number of marks transversal to the
boundary (i.e. such that the corresponding separating curve goes
inside the region) is even since at each of this marks the
orientation of the edges, along the boundary, is reverse, whereas
on the other crossings and vertices adjacent to the region the
orientation is preserved since the orientation is hyperbolic (see
fig.).

\centerline{\psbild{25}{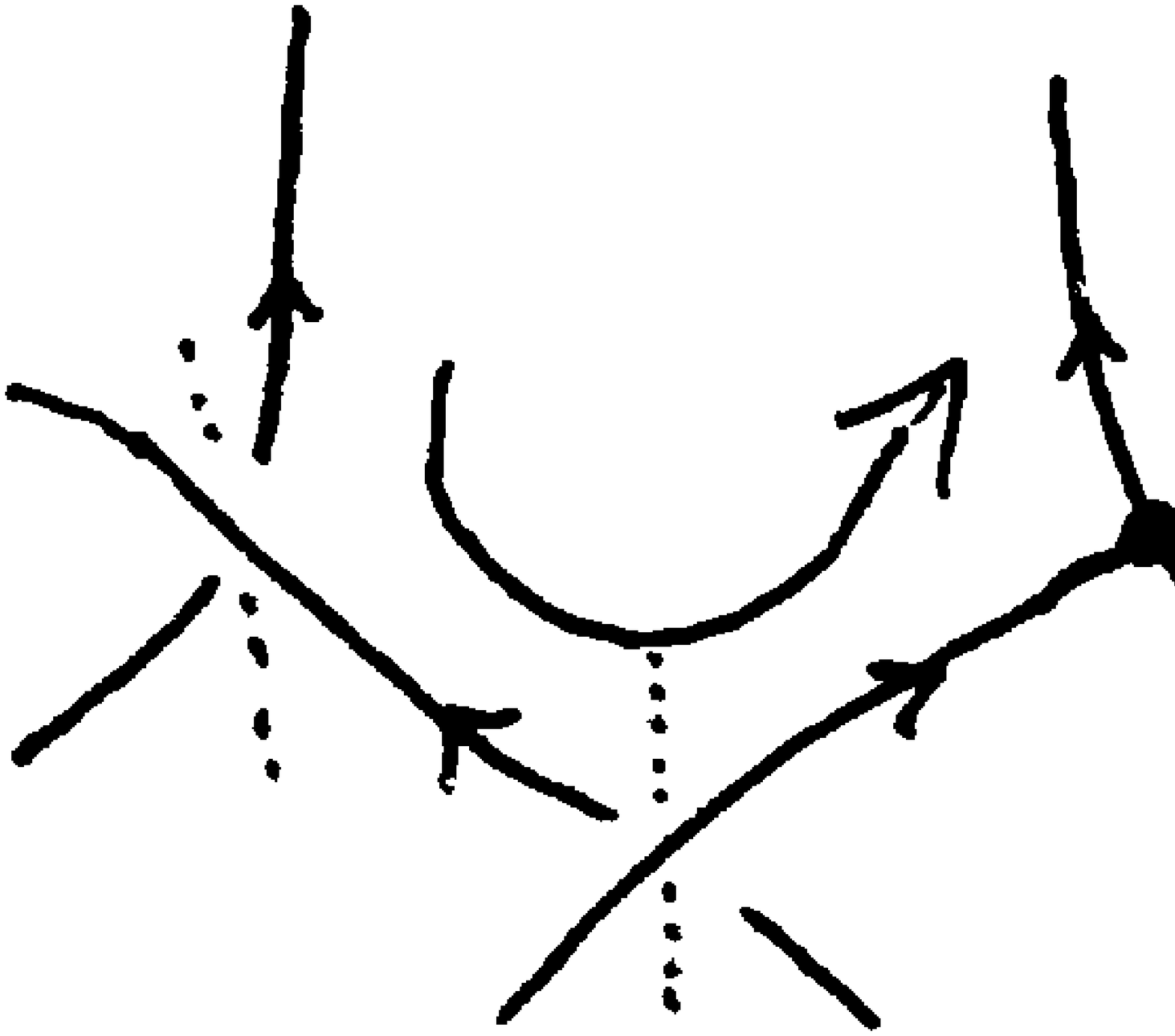}}

Thus, it is possible to join up such marks two by two using paths
without intersections.

$$ \hspace{1mm}\raisebox{-14mm}{\epsfysize28mm
\epsffile{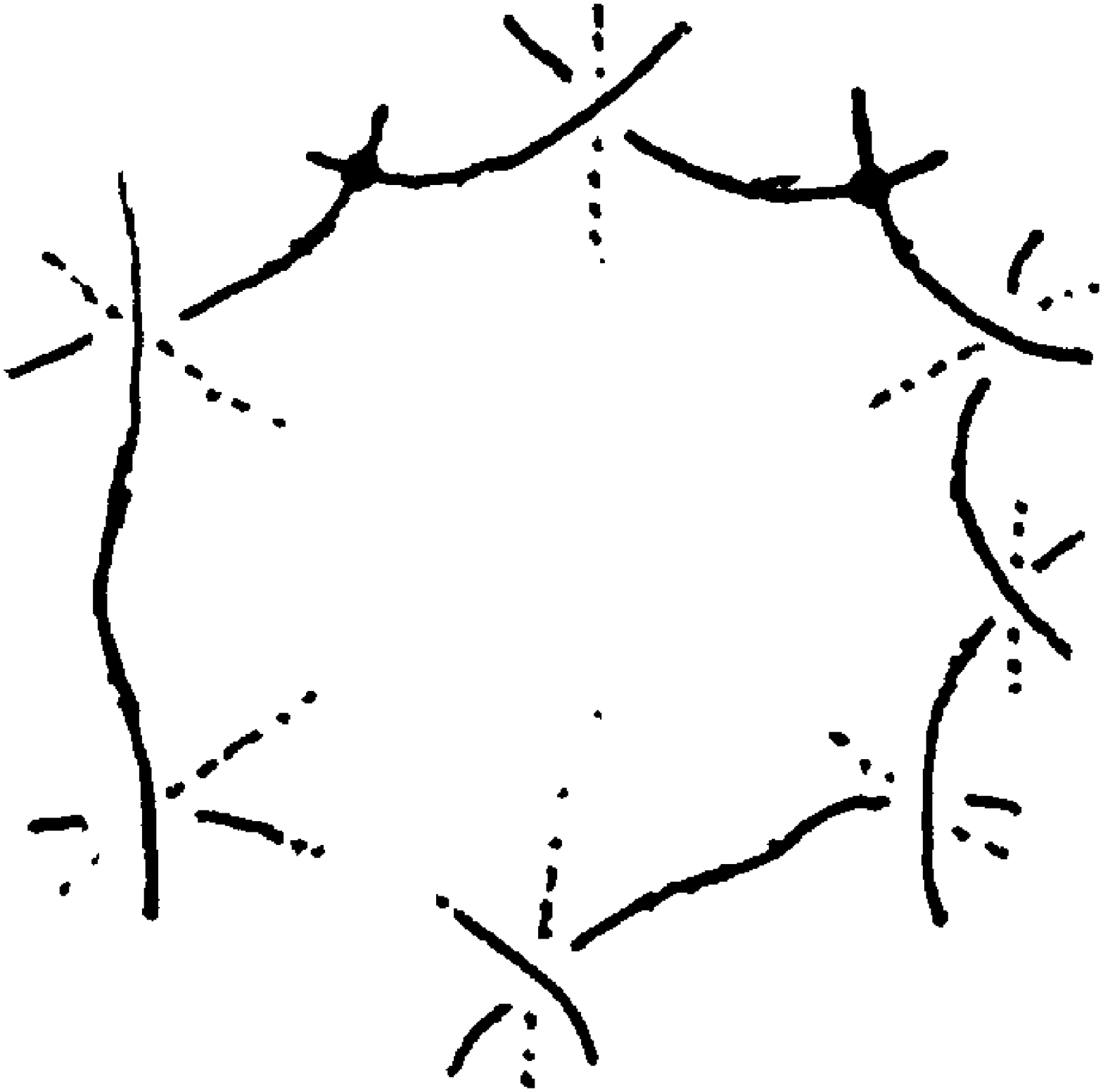}}\hspace{1mm} \longrightarrow
\raisebox{-14mm}{\epsfysize28mm \epsffile{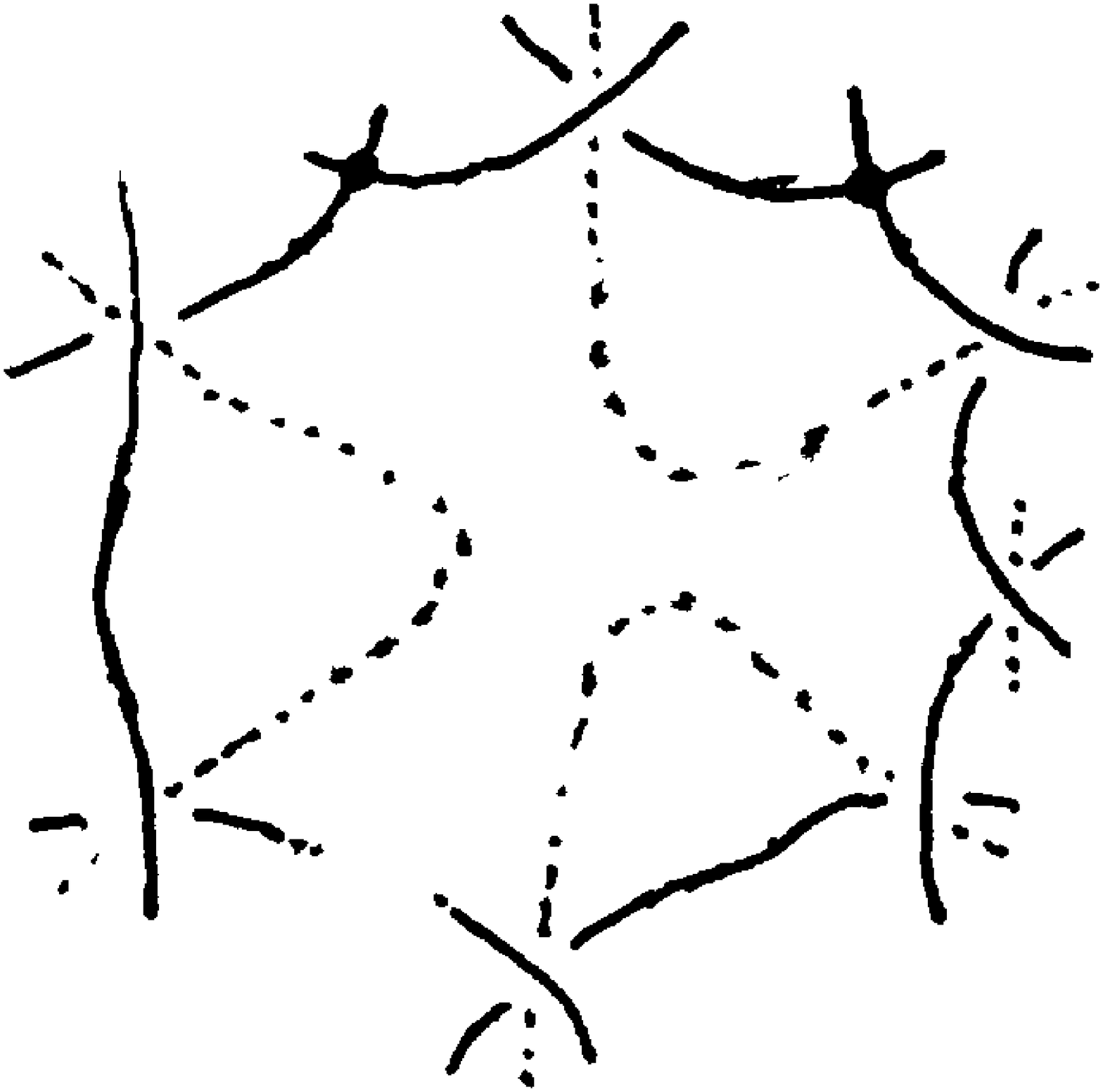}}\hspace{1mm}
$$

In the end we obtain a collection of separating curves for these
marks, so we get a partition. \TeXButton{End Proof}{\endproof}

\vspace{5mm}

For a partition, we say that the crossing \psdiag{2}{6}{a1} has
signature $A$ if its partition mark is vertical, and has signature
$A^{-1}$ if its partition mark is horizontal. The product of the
signatures of all crossings of the partition we call the {\it
signature} of the partition. By convention we fix the value $1$
for the signature of the null partition. For a graph diagram $G$
we denote by $\{G\}$ the sum of the signatures of all partitions
of the diagram. If $G$ has no partition (it is non separable) we
fix $\{G\}=0$.

Let $G$ be a graph diagram with a hyperbolic orientation on the
edges, and let $w(G)$ be the {\it writhe of G}, that is, the sum
of the signs of all crossings on this orientation. It is easy to
see that the signature of the partition corresponding to the given
orientation is $A^{w(D)}$.

A {\it diagram component} is a component corresponding to a
connected component of the planar graph obtained by changing
crossings into vertices. For example, a diagram of Hopf's link has
only one diagram component although Hopf's link has two connected
components.

We can also see that two hyperbolic orientations of the graph give
us the same partition if and only if, on each diagram component,
the orientations of the connected components of the graph are the
same or opposite\footnote{This operation changes neither the
writhe nor the hyperbolicity of the orientation.}. Moreover, we
have that, on each connected component of a graph, one hyperbolic
orientation, if it exists, it is determined by the orientation of
one edge in such a component.

Thus, calling ${\cal H}$ the set of all hyperbolic orientations of
the diagram $G$, and $G^{h}$ the diagram $G$ with a given
hyperbolic orientation $h$, we have $$\sum_{h \in {\cal
H}}A^{w(G^{h})}=2^{d}\{G\}$$ where $d$ is the number of diagram
components of $G$.

\begin{theorem}
The polynomial with variables $B=A^{-1}$ and $a=A$ for a graph
diagram $G$ is $$\frac{1}{2}(-A-A^{-1})^{v}\sum_{h \in {\cal
H}}A^{w(G^{h})}$$ where ${\cal H}$ is the set of all hyperbolic
orientations of the diagram $G$, $G^{h}$ is the diagram $G$ with a
given hyperbolic orientation $h$ and $v$ is the number of vertices
of the graph $G$.
\end{theorem}

\TeXButton{Proof}{\proof} We only need to check that the
polynomial $$[G]=\frac{1}{2}(-A-A^{-1})^{v}\sum_{h \in {\cal
H}}A^{w(G^{h})}$$ satisfies the skein relation and that, for
planar graphs, it is equal to $2^{c-1}(-A-A^{-1})^v$, where $c$ is
the number of connected components of $G$ and $v$ is the number of
vertices of $G$.

First of all, we check the skein relation (with $B=A^{-1}$ and
$a=A$).

Let [\psdiag{2}{6}{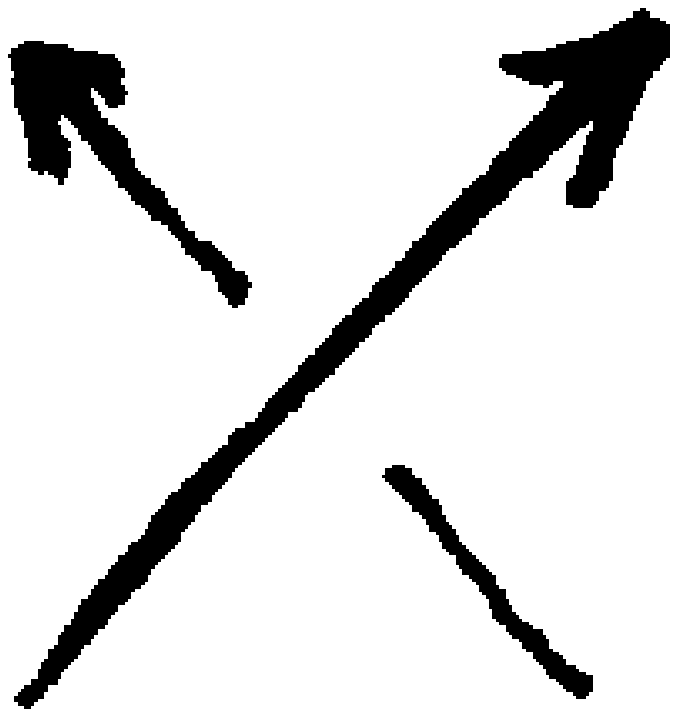}], ..., [\psdiag{2}{6}{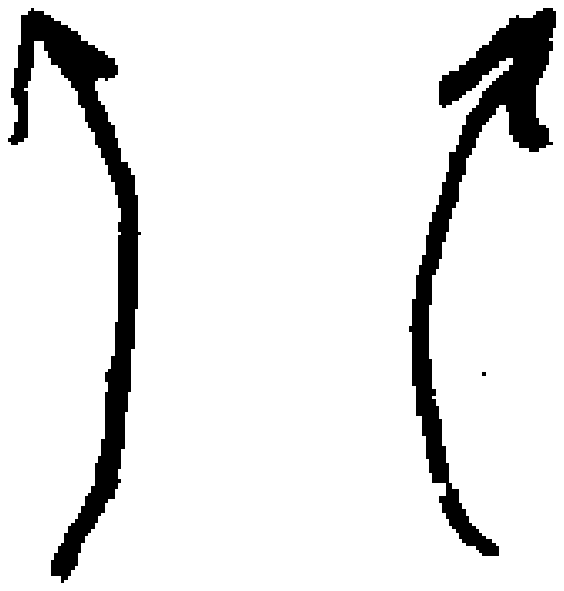}], ... be the
polynomials [\psdiag{2}{6}{a1}], [\psdiag{2}{6}{d1}], ...
considering only the hyperbolic orientation with the local form
\psdiag{2}{6}{w1}, ..., \psdiag{2}{6}{y1}, ... .Then we have the
following identities:

$$
\begin{array}{ccccccccccccc}
\lbrack
\psdiag{2}{6}{a1}]&=&[\psdiag{2}{6}{w1}]&+&[\psdiag{2}{6}{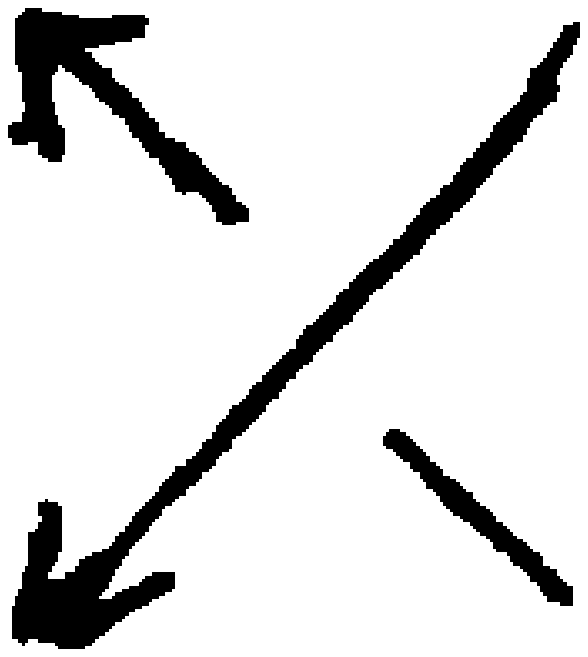}]&+&[\psdiag{2}{6}{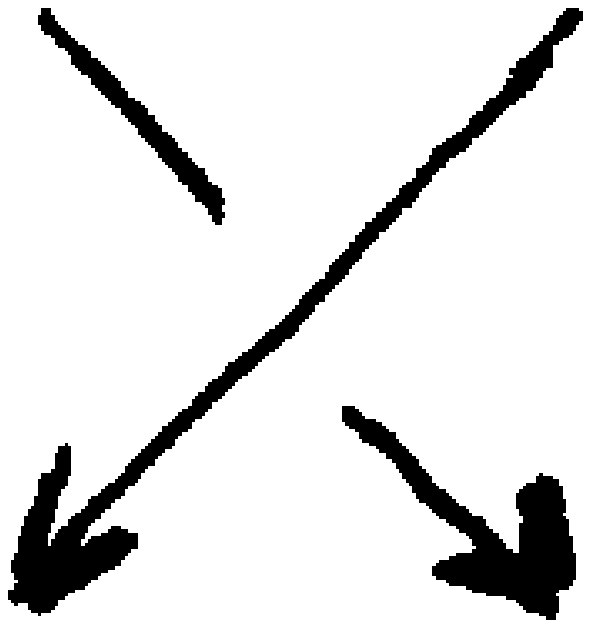}]&
+&[\psdiag{2}{6}{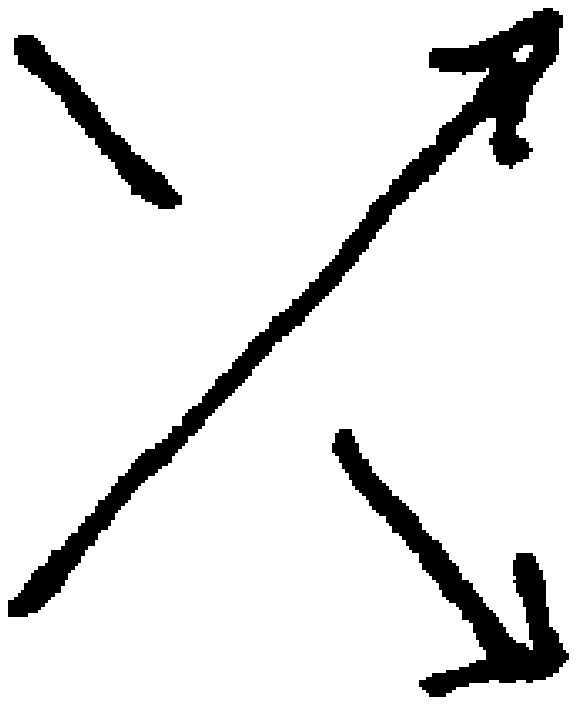}]&&&& \\ \lbrack
\psdiag{2}{6}{d1}]&=&[\psdiag{2}{6}{y1}]& &
&+&[\psdiag{2}{6}{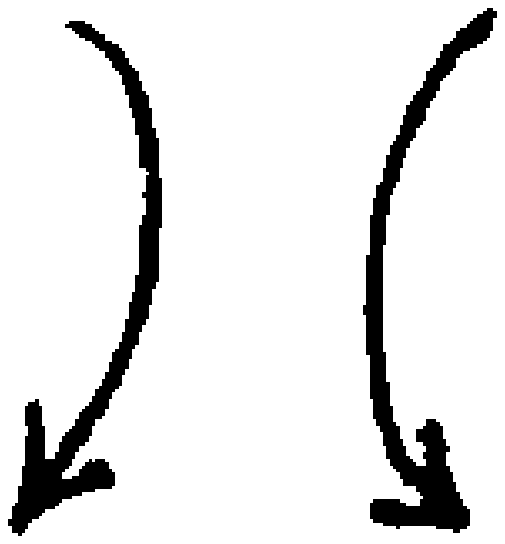}]& &
&+&[\psdiag{2}{6}{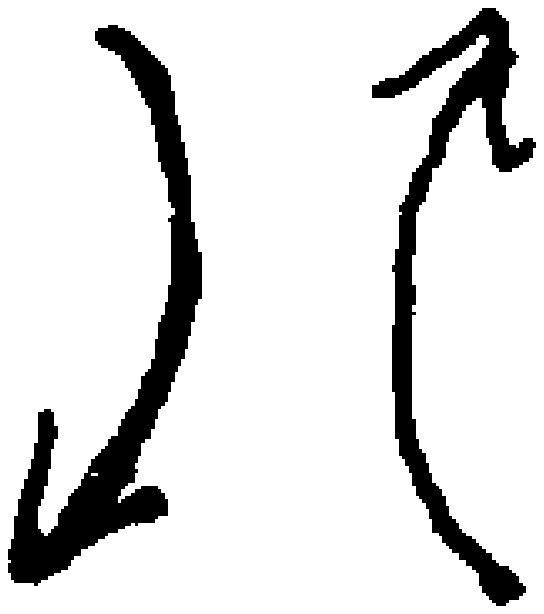}]&+&[\psdiag{2}{6}{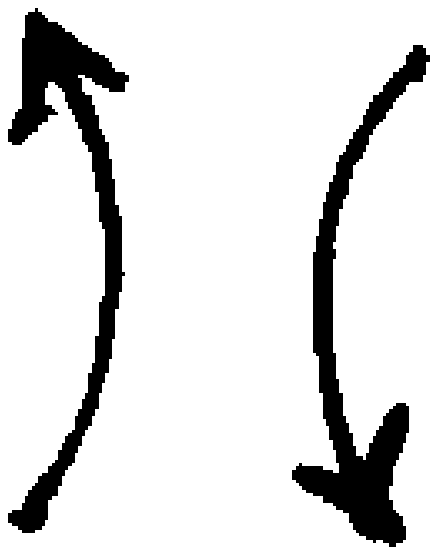}]
\\ \lbrack
\psdiag{2}{6}{e1}]&=& & &[\psdiag{2}{6}{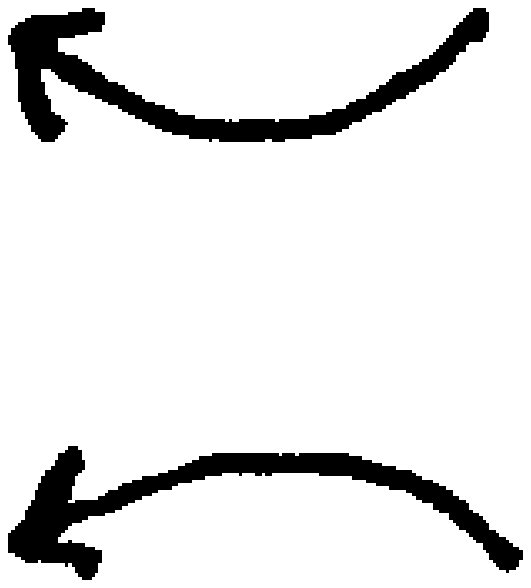}]& &
&+&[\psdiag{2}{6}{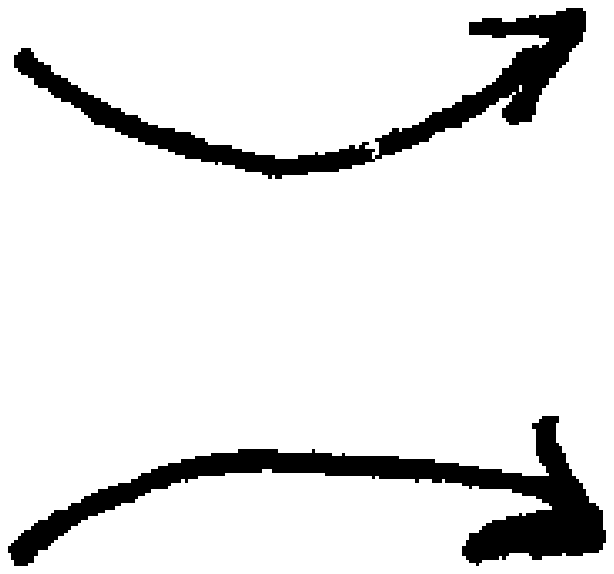}]&+&[\psdiag{2}{6}{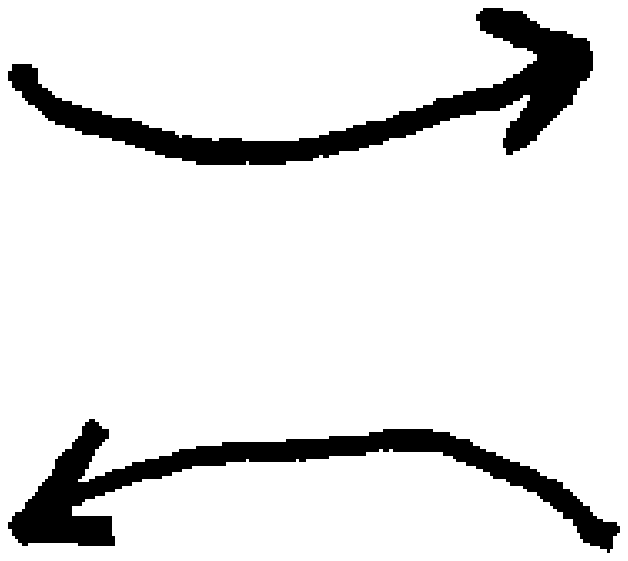}]&+&[\psdiag{2}{6}{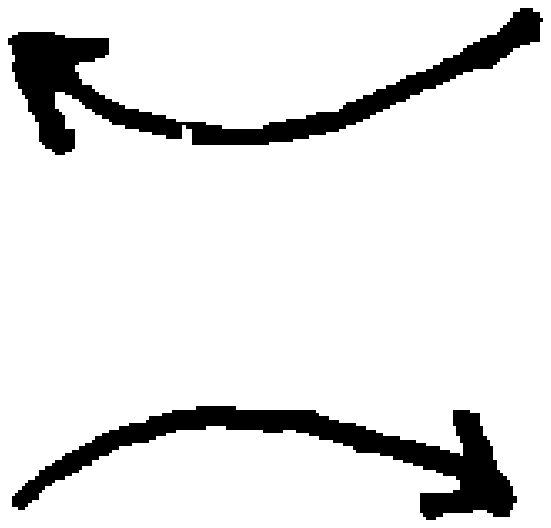}]
\\ \lbrack
\psdiag{2}{6}{c1}]&=& & & & & & & &
&[\psdiag{2}{6}{x8}]&+&[\psdiag{2}{6}{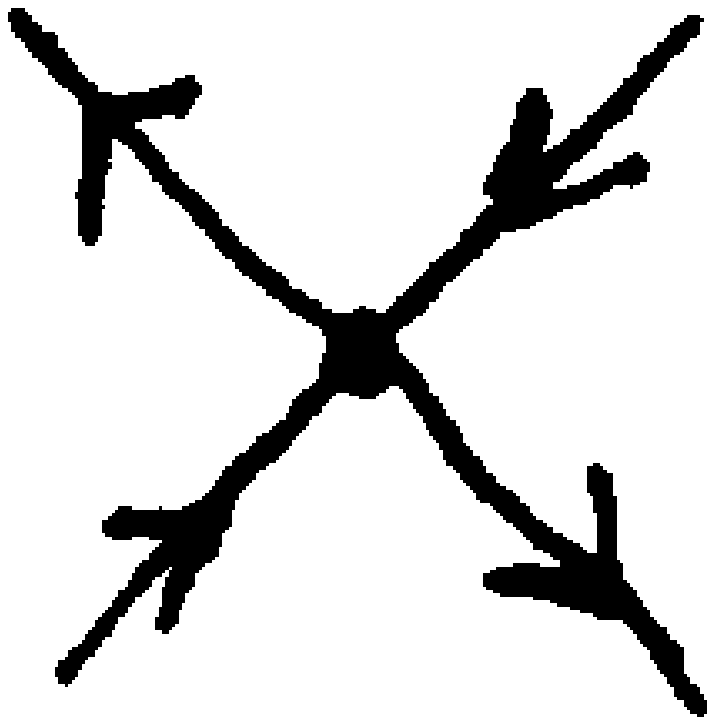}]
\end{array}
$$ and $$
\begin{array}{cccc}
\lbrack \psdiag{2}{6}{y4}]&=&[\psdiag{2}{6}{z4}]& \\ \lbrack
\psdiag{2}{6}{y3}]&=&[\psdiag{2}{6}{z3}]& \\ \lbrack
\psdiag{2}{6}{x8}]&=&-(A+A^{-1})&[\psdiag{2}{6}{y4}] \\ \lbrack
\psdiag{2}{6}{y8}]&=&-(A+A^{-1})&[\psdiag{2}{6}{y3}] \\ \lbrack
\psdiag{2}{6}{w1}]&=&A[\psdiag{2}{6}{y1}]& \\ \lbrack
\psdiag{2}{6}{w2}]&=&A^{-1}[\psdiag{2}{6}{z2}]& \\ \lbrack
\psdiag{2}{6}{w3}]&=&A[\psdiag{2}{6}{y2}]& \\ \lbrack
\psdiag{2}{6}{w4}]&=&A^{-1}[\psdiag{2}{6}{z1}]&
\end{array}
$$

Thus we take:

\begin{eqnarray*}
\lefteqn{A[\psdiag{2}{6}{d1}]+A^{-1}[\psdiag{2}{6}{e1}]+[\psdiag{2}{6}{c1}]=}\\
&=& A([\psdiag{2}{6}{y1}]+[\psdiag{2}{6}{y4}]+[\psdiag{2}{6}{y2}]
+[\psdiag{2}{6}{y3}])\\ &&
+A^{-1}([\psdiag{2}{6}{z4}]+[\psdiag{2}{6}{z2}]+[\psdiag{2}{6}{z3}]
+[\psdiag{2}{6}{z1}])\\ &&
+[\psdiag{2}{6}{x8}]+[\psdiag{2}{6}{y8}]=\\ &=&
A([\psdiag{2}{6}{y1}]+[\psdiag{2}{6}{y2}])+A^{-1}([\psdiag{2}{6}{z2}]+[\psdiag{2}{6}{z1}])+\\
&&
A[\psdiag{2}{6}{y4}]+A^{-1}[\psdiag{2}{6}{z4}]+[\psdiag{2}{6}{x8}]+A[\psdiag{2}{6}{y3}]
+A^{-1}[\psdiag{2}{6}{z3}]+[\psdiag{2}{6}{y8}]=\\ &=&
A([\psdiag{2}{6}{y1}]+[\psdiag{2}{6}{y2}])+A^{-1}([\psdiag{2}{6}{z2}]+[\psdiag{2}{6}{z1}])=\\
&=& [\psdiag{2}{6}{w1}]
+[\psdiag{2}{6}{w2}]+[\psdiag{2}{6}{w3}]+[\psdiag{2}{6}{w4}]= \\
&=& [\psdiag{2}{6}{a1}]
\end{eqnarray*}

Since, for a planar graph $G$, the writhe of any hyperbolic
orientation is zero, we have that $\sum_{h \in {\cal
H}}A^{w(G^{h})}$ is equal to the number of hyperbolic orientations
which is $2^{c}$ where $c$ is the number of connected components
of $G$. Thus the polynomial $[G]=\frac{1}{2}(-A-A^{-1})^{v}\sum_{h
\in {\cal H}}A^{w(G^{h})}$ is equal to $2^{c-1}(-A-A^{-1})^v$ for
a planar graph $G$. \TeXButton{End Proof}{\endproof}

\vspace{5mm}

Since any orientation on a link is hyperbolic, we have that all
links, including the knots, are separable.

In the case of the  knots, we have only two orientations and these
have the same writhe. Thus the polynomial for a knot $K$ is
$[K]=A^{w(K)}$ where $w(K)$ is the writhe of $K$. In this way, the
corrected polynomial $A^{-t(G)}[G]$ become trivial for knots.

Another result that is easy to check is the following.

\begin{corollary}
If the graph admits a knot-theoretic circuit that passes through
vertices an odd number of times then the graph is not separable.
\end{corollary}

This happens because a hyperbolic orientation reverses the
orientations of the edges along a knot-theoretic circuit, whenever
it passes through a vertex. This gives us a easy way to find non
separable graphs (though it is not general).

\section{Conjectures.}

Now we consider the polynomial $P$ that for a graph $G$ gives the
polynomial $$P(G)=\frac{A^{-t(G)}[G]}{2^{c-1}(-A-A^{-1})^v}.$$
This polynomial is invariant under isotopies and for planar graphs
takes the value 1. Moreover, this polynomial satisfies the
following identities:

\begin{description}
\item $P(\psdiag{2}{6}{k2})=P(\psdiag{2}{6}{g1})$;
\item $P(\psdiag{2}{6}{l2})=P(\psdiag{2}{6}{g1})$;
\item $P(\psdiag{2}{6}{m2})=P(\psdiag{2}{6}{d1})$, if
\psdiag{2}{6}{m2} and \psdiag{2}{6}{d1} admit hyperbolic
orientations and \psdiag{2}{6}{d1} has a hyperbolic orientation
with the local form \psdiag{2}{6}{y3} or \psdiag{2}{6}{y4};
\item $P(\psdiag{2}{6}{b3})=P(\psdiag{2}{6}{a3})$.
\end{description}

We are interested in studying the class of the graphs with the
polynomial $P$ equal to 1 (that is the class of graphs for which
the polynomials $P$ or $[\ ]$ does not detect non-planarity). This
class contains the planar graphs and the knots. Moreover, it
contains any graph that can be obtained from a planar graph or a
knot by the moves:

\begin{description}
\item $\psdiag{2}{6}{k2} \rightleftharpoons \psdiag{2}{6}{g1}$;
\item $\psdiag{2}{6}{l2} \rightleftharpoons \psdiag{2}{6}{g1}$;
\item $\psdiag{2}{6}{m2} \rightleftharpoons \psdiag{2}{6}{d1}$, if
\psdiag{2}{6}{m2} and \psdiag{2}{6}{d1} are hyperbolic and
\psdiag{2}{6}{d1} has a hyperbolic orientation with the local form
\psdiag{2}{6}{y3} or \psdiag{2}{6}{y4};
\item $\psdiag{2}{6}{b3} \rightleftharpoons \psdiag{2}{6}{a3}$.
\end{description}

Let ${\cal C}_1$ be the set of the graphs isotopic to planar
graphs, ${\cal C}_2$ the set of the graphs isotopic to graphs
obtained from planar graphs by these last moves, ${\cal C}_3$ the
set of the graphs isotopic to graphs obtained from
knots\footnote{including the unknot.} by these last moves, and
${\cal C}_4$ the set of the graphs with the polynomial $P$ equal
to 1.

We have that ${\cal C}_1 \subseteq {\cal C}_2$, ${\cal C}_2
\subseteq {\cal C}_4$ and ${\cal C}_3 \subseteq {\cal C}_4$.

We have also that ${\cal C}_2 \subseteq {\cal C}_3$ because, by
the lemma 2 of \cite{4}, any 4-valent planar graph can be
transformed into a disjoint union of circles (that are the knots
theoretic circuits of the graph) by the moves: $$\psdiag{2}{6}{l2}
\rightarrow \psdiag{2}{6}{g1}$$, $$\psdiag{2}{6}{m2} \rightarrow
\psdiag{2}{6}{d1}$$ and $$\psdiag{2}{6}{b3} \rightarrow
\psdiag{2}{6}{a3}$$ and by the move $$\psdiag{2}{6}{k2}
\rightarrow \psdiag{2}{6}{g1}$$ we obtain the unknot.

\begin{quote}
{\it Question: Are any two of the sets ${\cal C}_1$, ${\cal C}_2$,
${\cal C}_3$ and ${\cal C}_4$ the same?}
\end{quote}

We have that ${\cal C}_1  \neq {\cal C}_2$, as the following
example shows:

$$ G:\psdiag{14}{28}{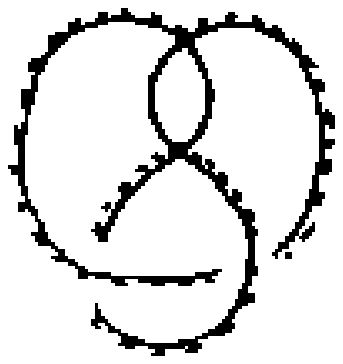} $$

\begin{example}
The graph in the picture is not planar, since it contains two
disjoint linked cycles (indicated in the picture), but it can be
obtained from a planar graph by the move $\psdiag{2}{6}{d1}
\rightarrow \psdiag{2}{6}{m2}$.
\end{example}

\begin{conjecture}
${\cal C}_2 \neq {\cal C}_3$.
\end{conjecture}

It is hard to believe that any knot can be obtained from a planar
graph by the moves:

\begin{description}
\item $\psdiag{2}{6}{k2} \rightleftharpoons \psdiag{2}{6}{g1}$;
\item $\psdiag{2}{6}{l2} \rightleftharpoons \psdiag{2}{6}{g1}$;
\item $\psdiag{2}{6}{m2} \rightleftharpoons \psdiag{2}{6}{d1}$, if
\psdiag{2}{6}{m2} and \psdiag{2}{6}{d1} are hyperbolic and
\psdiag{2}{6}{d1} has a hyperbolic orientation with the local form
\psdiag{2}{6}{y3} or \psdiag{2}{6}{y4};
\item $\psdiag{2}{6}{b3} \rightleftharpoons \psdiag{2}{6}{a3}$.
\end{description}

\begin{conjecture}
${\cal C}_3 \neq {\cal C}_4$.
\end{conjecture}

It is possible that the following example belongs to ${\cal C}_4$
but not to ${\cal C}_3$.

$$ G:\psdiag{14}{28}{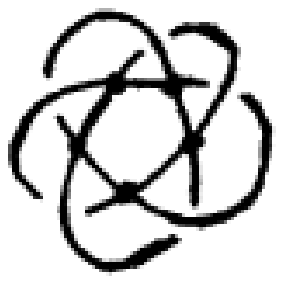} $$

\begin{example}
The graph in the picture is a representation of the complete graph
$K_5$ (thus is not planar), and has polynomial $P$ equal to 1.
\end{example}

{\bf Acknowledgment} -I wish to  thank Prof. Roger Picken who
encouraged me to write this paper and helped me to improve it.

This work was supported by the program {\em Programa Operacional
``Ci\^encia, Tecnologia, Inova\c c\~ao''} (POCTI) of the {\em
Funda\c c\~ao para a Ci\^encia e a Tecnologia} (FCT), cofinanced
by the European Community fund Feder.

\end{document}